\documentclass[11pt]{article}
\usepackage{amsmath}
\usepackage{amssymb}
\usepackage{graphicx}
\usepackage{amscd}
\usepackage{amsfonts}
\usepackage{bm}                                                
\usepackage{float}
\usepackage{latexsym}
\usepackage{lscape}
\usepackage{mathrsfs}
\usepackage{verbatim}

\newcommand{\ds}{\displaystyle}
\newtheorem{theorem}{Theorem}[section]

\newtheorem{lemma}{Lemma}[section]

\newtheorem{corollary}{Corollary}[section]

\newcommand{\mcI}{\ensuremath{\mathcal{I}}}

\newcommand{\mcS}{\ensuremath{\mathcal{S}}}

\newcommand{\mrI}{\ensuremath{\mathrm{I}}}

\def\qed{\hbox{\vrule width 6pt height 6pt depth 0pt}}

\title{Wellposedness of the two-sided variable coefficient Caputo flux fractional diffusion equation and error estimate of its spectral approximation} 

\author{
	Xiangcheng Zheng \thanks{Department of Mathematics, University of South Carolina, Columbia,
		South Carolina 29208, USA. email: {\tt xz3@math.sc.edu \& hwang@math.sc.edu}.} 
	\and
	V.~J.~Ervin\thanks{Department of Mathematical Sciences,
	  Clemson University, Clemson, South Carolina 29634-0975, USA.
	  email: {\tt vjervin@clemson.edu}. }
	\and 
	 Hong Wang $^*$ } 

\date{\today}

\begin{document}

\maketitle

\begin{abstract}
In this article a two-sided variable coefficient fractional diffusion equation (FDE) is investigated,
where the variable coefficient occurs outside of the fractional integral operator. Under a suitable transformation the variable
coefficient equation is transformed to a constant coefficient equation. Then, using
the spectral decomposition approach with Jacobi polynomials, we proved the wellposedness of the model and the 
regularity of its solution. A spectral approximation scheme is proposed and the accuracy of its approximation studied.
Two numerical experiments are presented to
demonstrate the derived error estimates.
\end{abstract}



\setcounter{equation}{0}
\setcounter{figure}{0}
\setcounter{table}{0}
\setcounter{theorem}{0}
\setcounter{lemma}{0}
\setcounter{corollary}{0}
\section{Introduction}

Fractional differential equations (FDEs) provide a competitive means to model challenging phenomena such as anomalously diffusive transport and long-range spatial interactions and memory effect. In the last couple of decades they have found increasingly more applications \cite{BenWhe,MetKla00} and have attracted extensive research on the development and numerical analysis of their numerical approximations. 

In their pioneer work Ervin and Roop \cite{ErvRoo05} proved the equivalence between the fractional derivative spaces and corresponding fractional Sobolev spaces and proved that the Galerkin weak formulation for the homogeneous Dirichlet boundary-value problem of a one-dimensional steady-state FDE with a constant diffusivity is coercive and bounded on the fractional Sobolev space $H^{\alpha/2}_0(0,1) \times H^{\alpha/2}_0(0,1)$, where $\alpha$ is the order of the fractional diffusion operator. These results enable them to apply Lax-Milgram theorem to conclude the wellposedness of the Galerkin weak formulation \cite{Cia}. In addition, they proved an optimal-order error estimate in the energy norm for the corresponding Galerkin finite element method under the assumption that the true solution of the FDE has full regularity. Furthermore, they used Nitsche lifting to prove an optimal-order error estimate in the $L^2$ norm for the finite element method under the assumption that the true solution to the dual problem has the full regularity for {\it each} right-hand side source term in the $L^2$. Subsequently, the framework in \cite{ErvRoo05} was extended by various researchers to prove error estimates for other numerical methods, again, under the assumption that the true solutions to these problems have full regularities. 

It was shown in \cite{WanYanZhu14,WanZha} that the solution to the FDE with a smooth right-hand side is actually not in the Sobolev space $W^{1 \, , \, 1/(2-\alpha)}$, and, in particular, not in the Sobolev space $H^1$ for $1 < \alpha < 3/2$. This shows that FDEs have significantly different mathematical behavior than their integer-order analogues \cite{Eva,GilTru}. Consequently, any Nitsche-lifting based proof of optimal-order $L^2$ error estimates of finite element methods in the literature does not hold since the full regularity assumption of the solution to the dual problem is invalid. Numerical experiments confirmed that spectral methods or high-order finite element methods for this problem failed to achieve high-order convergence rates in the $L_2$ norm \cite{WanYanZhu17,WanZha}. Jin et al. \cite{JinLaz} studied the regularity of the solution to a one-sided simplification of the problem and proved (suboptimal-order) error estimates of the corresponding Galerkin and Petrov-Galerkin finite element method in the Sobolev spaces only under the assumption of the smoothness of the data (not the solution). 

It turns out that spectral methods are particularly well suited for the accurate approximation of FDEs for the following reasons: (i) They present a clean analytical expression of the true solution to FDEs, which have been fully explored in \cite{ErvHeuRoo,MaoGeo} in analyzing the structure and regularity of the true solutions. (ii) Fractional differentiation of many spectral polynomials can be carried out analytically \cite{WanZha}, in contrast to finite element methods in which they have to be calculated numerically that are sometimes a headache \cite{WanYanZhu17}. (iii) As FDEs are nonlocal operators the appealing property of a sparse coefficient matrix, which arises in a finite element, finite difference, or finite volume approximation of a usual differential equation, is lost. In contrast, the stiffness matrices of spectral methods are often diagonal (at least for constant-coefficient FDEs). Mao et al. \cite{MaoCheShe} analyzed the regularity of the solution to a symmetric case of the FDE and developed corresponding spectral methods. The solution structure to the general case was resolved completely in \cite{ErvHeuRoo}, in which a spectral method utilizing the weighted Jacobi polynomial was studied and a priori error estimates derived. The two-sided FDE with constant coefficient and Riemann-Liouville fractional derivative was investigated in \cite{MaoGeo}, by employing a Petrov-Galerkin projection in a properly weighted Sobolev space using two-sided Jacobi polyfractonomials as test and trial functions. 

The variable diffusivity $K$ presents another bottleneck of FDEs. It was shown in \cite{WanYan} that the Galerkin weak formulation may lose its coercivity for a smooth $K(x)$ with positive lower and upper bounds. In fact, its Galerkin finite element approximation might diverge \cite{WanYanZhu17}. A Petrov-Galerkin weak formulation was proved to be wellposed on $H^{\alpha-1}_0 \times H^1_0$ for $3/2 < \alpha < 2$ for a one-sided FDE with variable $K$ \cite{WanYan}. A Petrov-Galerkin finite element method was developed and analyzed subsequently \cite{WanYanZhu15}. An indirect Legendre spectral Galerkin method  and a finite element method  were developed for the FDE in \cite{WanZha} and \cite{WanYanZhu17}, respectively, in which the solution was expressed as a fractional derivative of the solution to a second-order differential equations. Consequently, high-order convergence rates of numerical approximations were proved and numerically observed in the $L_2$ norm only under regularity assumptions of coefficients and right-hand side (the true solution is not smooth in fact). In \cite{LiCheWan}, with the introduction of an auxiliary variable, a mixed approximation was developed and analyzed for an  FDE in which the variable diffusivity $K$ appears inside the fractional integral operator, and error estimates was derived. In \cite{MaoShe}, a spectral Galerkin method for a different variable coefficient FDE was analyzed, in which the outside and inside fractional derivatives are chosen carefully so that the corresponding Galerkin weak formulation are self-adjoint and coercive. Additionally optimal error estimates were derived under suitable smoothness assumption on the solution. 
Recently a FDE for a general $K(x)$ insider the fractional differential operator was studied in \cite{ZheErvWan}. The model was shown to be
wellposed, and a spectral approximation method studied.

In this paper we investigate the FDE where the variable diffusivity $K(x)$ appears outside the fractional integral operator. The approach
used in \cite{ZheErvWan} is not suitable for the solution of the FDE with $K(x)$ outside the fractional integral operator. Using a different
approach than in \cite{ZheErvWan}, we are able to establish wellposedness of the model. Additionally, a spectral approximation scheme
is proposed and analyzed.

This paper is organized as follows. In Section 2 we present the formulation of the model and introduce notations and key lemmas used in the analysis. The wellposedness of the model and the regularity of its solution are studied in Section 3. The spectral approximation method is formulated and a detailed analysis of its convergence is given in Section 4. Two numerical experiments are presented in Section 5 whose results demonstrate the sharpness of the derived error estimates.

\setcounter{equation}{0}
\setcounter{figure}{0}
\setcounter{table}{0}
\setcounter{theorem}{0}
\setcounter{lemma}{0}
\setcounter{corollary}{0}

\section{Model problem and preliminaries}

In this paper we consider the following homogeneous Dirichlet boundary-value problem of a two-sided Caputo flux FDE, which is obtained by incorporating a fractional Fick's law into a conventional local mass balance law \cite{CasCar,ErvRoo05,ZhaBen}. 
Fractional diffusion equation (FDE) of order $1 < \alpha < 2$ with variable fractional diffusivity \cite{CasCar,ErvRoo05,ZhaBen} 
\begin{equation}\label{ModelV}\begin{array}{c}
 -D \big [K(x) \big (r \,{}_0I_x^{2-\alpha} +  (1-r) \,  {}_xI_1^{2-\alpha} \big ) \, Du(x) \big ] = f(x), ~~x \in (0,1), \\[0.05in]
u(0) = u(1) = 0, \qquad 1 < \alpha < 2.
\end{array}\end{equation}
Here $D$ is the first-order differential operator, $K(x)$ is the fractional diffusivity with $0 <  K_{m} \leq K(x) \leq  K_{M} < \infty$, $0\leq r\leq 1$ indicates the relative weight of forward versus backward transition probability and $f(x)$ the source or sink term. The left and right fractional integrals of order $0 < \sigma < 1$ are defined as \cite{Pod}
\begin{equation*}
{}_0I_x^{\sigma}w(x) := \frac{1}{\Gamma(\sigma)} \int_0^x w(s)(x-s)^{\sigma-1}ds, \quad {}_xI^{\sigma}_1w(x) := \frac{1}{\Gamma(\sigma)} \int_x^1 w(s)(s-x)^{\sigma-1} \, ds \, ,
\end{equation*}
where $\Gamma(\cdot)$ is the Gamma function. 

\begin{lemma} \label{lem:EquivForm}
For $f(x) \in L_{1}(0,1)$,
$u$ is a solution to problem (\ref{ModelV}) if and only if $u$ is a solution to the following problem
\begin{equation}\label{EquivForm}\begin{array}{c}
\ds L_r^{\alpha}u(x) := -\big(r{}_0I_x^{2-\alpha}+(1-r){}_xI_1^{2-\alpha}\big)Du(x) = f_1(x) - A f_2(x), ~~x \in (0,1),\\[0.05in]
~~ u(0) = u(1) = 0,
\end{array}\end{equation}
where 
\begin{equation}\label{EquivForm:e1}\begin{array}{rl}
A & \ds := K(0)\big[\big(r{}_0I_x^{2-\alpha} + (1-r){}_xI_1^{2-\alpha}\big)Du(x)\big]\big|_{x=0}, \\[0.1in]
f_1(x) &\ds := \frac1{K(x)} \int_0^x f(y)dy, \quad f_2(x) := \frac{1}{K(x)}. 
\end{array}\end{equation}
\end{lemma}
{\bf Proof.} Note that integration of the FDE (\ref{ModelV}) with simple algebraic manipulation yields (\ref{EquivForm}) and the procedure is reversible. \qed

Based on Lemma \ref{lem:EquivForm}, 
in the subsequent sections we concentrate on analyzing problem (\ref{EquivForm}). In the rest of this section we introduce some notion and known results in the literature to be used in subsequent sections. 

The notation $y_{n} \sim n^{p}$ means that there exist constants $0 < C_1 < C_2 < \infty$ such that $C_1 n^{p} \le y_{n} \le C_2 n^{p}$.

We use $C$ to denote a generic positive constant, whose actual value may change from line to line.

Let $\mrI  :=   (0 , 1)$, $\mathbb{N}_{0}  := \mathbb{N} \cup {0}$,
and $W^r_p(\mrI)$ for $r\geq 0,~1\leq p\leq \infty$ denote the standard Sobolev spaces with the corresponding norms $\|\cdot\|_{W^r_p}:=\|\cdot\|_{W^r_p(\mrI)}$. For $\omega(x) > 0, ~x \in \mrI$, we introduce the weighted $L^2_{\omega}$
inner product and associated norm defined by
\[
\left( f \, , \, g \right)_{\omega} \, := \, \int_{0}^{1} \omega(x) \, f(x) \, g(x) \, dx \, , \quad \quad
\| f \|_{\omega} \, := \, \left( f \, , \, f \right)_{\omega}^{1/2} \, .
\]
Correspondingly, we let $L_{\omega}^{2}(\mrI)$ denote the weighted $L^{2}$ space
$$\ds L_{\omega}^{2}(\mrI) \, := \, \bigg \{ f(x) \, : \, \| f\|_{\omega}^2  \ < \ \infty \bigg \}.$$

In addition, let $\omega^{(\alpha,\beta)}$ be a weighting function defined on $\mrI$ and indexed by $\alpha$ and $\beta$. For any $m \in \mathbb{N}$, we introduce the following weighted Sobolev spaces \cite{GuoWan,SheTan}
\begin{equation}\label{Hw}
\ds H^m_{\omega^{(\alpha,\beta)} }(\mrI) := \bigg \{ v \, : \| v \|_{m,\omega^{(\alpha,\beta)}}^2 := \sum_{j=0}^m \big | v \big |_{j,\omega^{(\alpha,\beta)}}^2 = \sum_{j=0}^m \big \| D^j v \big \|_{\omega^{(\alpha+j,\beta+j)}}^2 < \infty \bigg \}.
\end{equation}
For $t \in \mathbb{R}^{+}$ definition \eqref{Hw} is extended to $H^{t}_{\omega^{(\alpha ,\beta)}}(\mcI)$ by interpolation
\cite{BerLof}.

Jacobi polynomial play a key role in the approximation schemes. We briefly review their definition and
properties central to the method \cite{AbrSte, Sze}. 
  
\underline{Definition}: $ P_{n}^{(\alpha , \beta)}(x) \ := \ 
\sum_{m = 0}^{n} \, p_{n , m} \, (x - 1)^{(n - m)} (x + 1)^{m}$, where
\begin{equation}
       p_{n , m} \ := \ \frac{1}{2^{n}} \, \left( \begin{array}{c}
                                                              n + \alpha \\
                                                              m  \end{array} \right) \,
                                                    \left( \begin{array}{c}
                                                              n + \beta \\
                                                              n - m  \end{array} \right) \, .
\label{spm21}
\end{equation}
\underline{Orthogonality}:    
\begin{align}
 & \int_{-1}^{1} (1 - x)^{\alpha} (1 + x)^{\beta} \, P_{j}^{(\alpha , \beta)}(x) \, P_{k}^{(\alpha , \beta)}(x)  \, dx 
 \ = \
   \left\{ \begin{array}{ll} 
   0 , & k \ne j  \\
   |\| P_{j}^{(\alpha , \beta)} |\|^{2}
   \, , & k = j  
    \end{array} \right.  \, ,  \nonumber \\
& \quad \quad \mbox{where } \  \ |\| P_{j}^{(\alpha , \beta)} |\| \ = \
 \bigg ( \frac{2^{(\alpha + \beta + 1)}}{(2j \, + \, \alpha \, + \, \beta \, + 1)} 
   \frac{\Gamma(j + \alpha + 1) \, \Gamma(j + \beta + 1)}{\Gamma(j + 1) \, \Gamma(j + \alpha + \beta + 1)}
   \bigg )^{1/2} \, .
  \label{spm22}
\end{align}         

In order to transform the domain of the family of Jacobi polynomials to $[0 , 1]$, let $x \rightarrow 2t - 1$ and 
introduce $G_{n}^{\alpha , \beta}(t) \, = \, P_{n}^{\alpha , \beta}( x(t) )$. From \eqref{spm22},
\begin{align}
 \int_{-1}^{1} (1 - x)^{\alpha} (1 + x)^{\beta} \, P_{j}^{(\alpha , \beta)}(x) \, P_{k}^{(\alpha , \beta)}(x)  \, dx 
 &= \
 \int_{t = 0}^{1} 2^{\alpha} \, (1 - t)^{\alpha} \, 2^{\beta} \, t^{\beta} \, P_{j}^{(\alpha , \beta)}(2t - 1) \, P_{k}^{(\alpha , \beta)}(2t - 1)  \, 2 \,  dt
 \nonumber \\
  &= \
2^{\alpha + \beta + 1} \int_{t = 0}^{1}   (1 - t)^{\alpha}  \, t^{\beta} \, G_{j}^{(\alpha , \beta)}(t) \, G_{k}^{(\alpha , \beta)}(t)  \,  dt
\nonumber \\
&= \
   \left\{ \begin{array}{ll} 
   0 , & k \ne j \, , \\
  2^{\alpha + \beta + 1} \, |\| G_{j}^{(\alpha , \beta)} |\|^{2}
   \, , & k = j  
   \, . \end{array} \right.    \nonumber \\
 \quad  \mbox{where } \  \  \big | \big \| G_{j}^{(\alpha , \beta)} \big | \big \| = \
 \bigg(& \frac{1}{(2j \, + \, \alpha \, + \, \beta \, + 1)} 
   \frac{\Gamma(j + \alpha + 1) \, \Gamma(j + \beta + 1)}{\Gamma(j + 1) \, \Gamma(j + \alpha + \beta + 1)}
   \bigg)^{1/2}  .  \label{spm22g} 
\end{align}                                                    

\begin{equation}
\mbox{Note that } \quad  |\| G_{j}^{(\alpha , \beta)} |\| \ = \ |\| G_{j}^{(\beta , \alpha)} |\| \, .
\label{nmeqG}
\end{equation}

 From \cite[equation (2.19)]{MaoCheShe} we have that
\begin{equation}
   \frac{d^{k}}{dx^{k}} P_{n}^{(\alpha , \beta)}(x) \ = \ 
   \frac{\Gamma(n + k + \alpha + \beta + 1)}{2^{k} \, \Gamma(n + \alpha + \beta + 1)} P_{n - k}^{(\alpha + k \, , \, \beta + k)}(x) \, .
   \label{derP}
\end{equation}   
Hence,
\begin{align}
\frac{d^{k}}{dt^{k}} G_{n}^{(\alpha , \beta)}(t) 
  &= \ \frac{\Gamma(n + k + \alpha + \beta + 1)}{  \Gamma(n + \alpha + \beta + 1)} 
  G_{n - k}^{(\alpha + k \, , \, \beta + k)}(t)  \, .  \label{eqC4}
\end{align}   

Also, from \cite[equation (2.15)]{MaoCheShe}, 
\begin{equation}
\begin{array}{l}
\ds\frac{d^{k}}{dx^{k}} \left\{ (1 - x)^{\alpha + k} \, (1 + x)^{\beta + k} \, P_{n - k}^{(\alpha + k \, , \, \beta + k)}(x) \right\}
\ \\[0.1in]
\ds\quad\quad\quad\quad\quad\quad= \ 
\frac{(-1)^{k} \, 2^{k} \, n!}{(n - k)!} \, (1 - x)^{\alpha} \, (1 + x)^{\beta} \, P_{n}^{(\alpha \, , \, \beta)}(x) \, , \
n \ge k \ge 0 \, ,
\label{eqB0}
\end{array}
\end{equation}
from which it follows that
\begin{equation}
 \frac{d^{k}}{dt^{k}} \left\{  \ (1 \, - \, t)^{\alpha + k} \,  t^{\beta + k} \,
 G_{n - k}^{(\alpha + k \, , \, \beta + k)}(t) \right\} 
 \ = \ 
 \frac{(-1)^{k} \,  n!}{(n - k)!} \,  (1 \, - \, t)^{\alpha} \,  t^{\beta} \,
 G_{n}^{(\alpha \, , \, \beta)}(t) \, . 
  \label{eqC2}
 \end{equation}
 
From \cite{ZheErvWan} we have
\begin{equation}\label{GnNormRelation}
\frac{1}{2} \, \le \, \frac{  |\| G_{j}^{(\alpha - \beta \, , \, \beta)}  |\|^{2} }{  |\| G_{j + 1}^{(\beta - 1 \, , \, \alpha - \beta - 1)} |\|^{2}} 
	\, = \, \frac{j + 1}{j + \alpha} \, \le \, 1, \, \quad j \ge 0.
\end{equation}

Let $\mcS_{N}$ denote the space of polynomials of degree $\le N$. We define the weighted $L^2$ orthogonal projection
$P_N^{(\alpha,\beta)} : \, L^{2}_{\omega^{(\alpha ,\beta)}}(0,1) \rightarrow \mcS_{N}$ 
\begin{equation}\label{Proj}
\big ( v \, - \,  P_N^{(\alpha,\beta)}v \ , \ \phi_N \big)_{\omega^{(\alpha ,\beta)}} \ = \ 0 \, , \ \ \forall \phi_N \in \mcS_{N}.
\end{equation}
\begin{lemma}\label{lem:Approx} \cite[Theorem 2.1]{GuoWan}
For any $v \in H^m_{\omega^{(\alpha ,\beta)}}(0,1)$, $m \in \mathbb{N}$, and $0 \le \mu \le m$, there exists a positive
constant $C$, which is independent of $N, \, \alpha$ or $\beta$, such that
\begin{equation}\label{Approx}
\big \| v  -  P_N^{(\alpha,\beta)} v \|_{\mu , \omega^{(\alpha ,\beta)}} \ \le \ C \, 
\big ( N \, (N + \alpha+\beta) \big )^{\frac{\mu - m}{2}} \, | v |_{m,\omega^{(\alpha ,\beta)}}.
\end{equation}
\end{lemma}
\mbox{ } \hfill \qed  

The Jacobi polynomials have the following useful property.
\begin{lemma}\label{lem:beta} \cite{MaoGeo} 
For $1 < \alpha < 2$ and $0 \le r \le 1$, let $\beta$, with $\alpha - 1 \le \beta \le 1$, be determined by
\begin{equation}\label{beta}
r \ = \ \frac{\sin( \pi \, \beta)}{\sin( \pi ( \alpha - \beta)) \, + \,  \sin( \pi \, \beta)}.  
\end{equation}
Then, for $n \, = \, 0, 1, 2, \ldots$
\begin{equation}\label{SturmLiouville}\begin{array}{c}
\ds \big(r\, {}_0I_x^{2-\alpha} + (1-r) \, {}_xI_1^{2-\alpha} \big) D\big(\omega^{(\alpha - \beta,\beta)}
\, G_{n}^{(\alpha - \beta \, , \, \beta)}(x) \big) = \ \lambda_{n}  \, G_{n + 1}^{(\beta - 1 \, , \, \alpha - \beta - 1)}(x), \\[0.1in]
\ds 
\mbox{where } \ \ \ \lambda_{n} = \ \frac{\sin(\pi \alpha)}{\sin(\pi (\alpha - \beta)) \ + \ \sin(\pi \beta)} \frac{\Gamma(n + \alpha)}{n !} \, .
\end{array}\end{equation}
\end{lemma}
\mbox{ } \hfill \qed 

By Stirling's formula, for $\mu > 0$, 
\begin{equation}
 \lim_{n \rightarrow \infty} \, \frac{\Gamma(n + \mu)}{\Gamma(n) \, n^{\mu}}
\ = \ 1 \, , \quad \mbox{hence }  \ \  \lambda_{n} \sim   (n + 1)^{\alpha - 1}.
\label{strfeq}
\end{equation}

\setcounter{equation}{0}
\setcounter{figure}{0}
\setcounter{table}{0}
\setcounter{theorem}{0}
\setcounter{lemma}{0}
\setcounter{corollary}{0}

\section{Wellposedness of \eqref{EquivForm} and the regularity of its solution}
\label{secspax}

In this section we study the wellposedness of \eqref{EquivForm} and the regularity of its solution.

\subsection{Existence and uniqueness of the solution to (\ref{EquivForm})}

\begin{lemma} \label{lem:f1f2}
For $1 < \alpha < 2$ and $0 \le r \le 1$, if $f \in L^{2}_{\omega^{(\beta , \alpha - \beta)}}(\mrI)$, then $\int_{0}^{x} f(s)ds, \, f_1, \, f_2 \in L^{2}_{\omega^{(\beta - 1,\alpha - \beta - 1)}}(\mrI)$, with the stability estimates
\begin{equation}\label{f1f2}
\| f_1 \|_{\omega^{(\beta-1, \alpha-\beta-1)}} \le C_{1}(\alpha,\beta, K_{m}) \| f \|_{\omega^{(\beta \, , \, \alpha-\beta)}}, \qquad \| f_2 \|_{\omega^{(\beta-1, \alpha-\beta-1)}} \le C_{2}(\alpha,\beta, K_{m}),
\end{equation}
for constants $C_{1}$ and $C_{2} > 0$.
\end{lemma}

\textbf{Proof}: 
To establish that $\int_{0}^{x} f(s)ds \,  \in L^{2}_{\omega^{(\beta - 1,\alpha - \beta - 1)}}(\mrI)$, consider
\begin{align}
&\int_0^1 \omega^{(\beta-1 \, , \,  \alpha-\beta-1)}(x) \Big (\int_0^x f(s)ds \Big)^2dx   \nonumber \\
& \quad \le \int_0^1\omega^{(\beta-1 \, , \, \alpha-\beta-1)}(x) \, \left( \int_0^x\omega^{(-\beta \, , \, -(\alpha-\beta))}(s)ds 
 \int_0^x\omega^{(\beta,\alpha-\beta)}f(s)^2 ds \right) dx  \nonumber \\
& \quad \le  \int_0^1\omega^{(\beta-1,\alpha-\beta-1)}(x) dx \,  \int_0^1\omega^{(-\beta \, , \, -(\alpha-\beta))}(s)ds \, 
 \int_0^1\omega^{(\beta,\alpha-\beta)}f(s)^2 ds   \nonumber \\
& \quad = \ B(\alpha-\beta \, , \, \beta) \, B(1-(\alpha-\beta) , \, 1-\beta) \, \| f \|_{\omega^{(\beta \, , \, \alpha-\beta)}}^{2} \, 
< \, \infty \, ,   \label{fbdn2}
\end{align} 
where $B(\cdot , \cdot)$ denotes the Beta function. Hence 
$\int_{0}^{x} f(s)ds \,  \in L^{2}_{\omega^{(\beta - 1,\alpha - \beta - 1)}}(\mrI)$.

For $f_{1}(x) \ = \ \frac{1}{K(x)} \int_{0}^{x} f(s)ds$, 
\begin{align*}
\| f_{1} \|_{\omega^{(\beta - 1,\alpha - \beta - 1)}} &= \
\|   \frac{1}{K(x)} \int_{0}^{x} f(s)ds   \|_{\omega^{(\beta - 1,\alpha - \beta - 1)}}
\ \le \ \frac{1}{K_{m}} \left( \int_0^1 \omega^{(\beta-1 \, , \,  \alpha-\beta-1)}(x) \Big (\int_0^x f(s)ds \Big)^2dx \right)^{1/2}
 \\
&\le \ \underbrace{ \frac{1}{K_{m}} \left( B(\alpha-\beta \, , \, \beta) \, B(1-(\alpha-\beta) , \, 1-\beta) \right)^{1/2}}%
_{ := \, C_{1}(\alpha,\beta, K_{m})} 
\| f \|_{\omega^{(\beta \, , \, \alpha-\beta)}} \, , \ \ \mbox{ using \eqref{fbdn2}} \, .
\end{align*}

Finally, for $f_{2}(x) \ = \ \frac{1}{K(x)}$,
\begin{align*}
\| f_{2} \|_{\omega^{(\beta - 1,\alpha - \beta - 1)}} &= \
\|   \frac{1}{K(x)}   \|_{\omega^{(\beta - 1,\alpha - \beta - 1)}}
\ \le \ \frac{1}{ K_{m}} \left( \int_0^1 \omega^{(\beta-1 \, , \,  \alpha-\beta-1)}(x) dx \right)^{1/2}  \\
&= \ \frac{1}{ K_{m}} \left( B(\alpha-\beta \, , \, \beta)  \right)^{1/2} \, := \, C_{2}(\alpha,\beta, K_{m}) \, .
\end{align*}
\mbox{ } \hfill \qed

\begin{theorem}\label{thm:Wellpose}
For $f \in L^{2}_{\omega^{(\beta , \alpha - \beta)}}(\mrI)$, problem (\ref{EquivForm}) has a solution 
$u \in L^{2}_{\omega^{(-(\alpha - \beta) ,  -\beta)}}(\mrI)$ in which the constant $A$ can be expressed in 
terms of the given data $f$ and $K$ as  
\begin{equation}\label{A}
A = \frac{\big (f_1,1 \big)_{\omega^{(\beta - 1,\alpha - \beta - 1)}} }{\big (K^{-1},1 \big)_{\omega^{(\beta - 1,\alpha - \beta - 1)}}} = \frac{\ds \bigg (\frac1{K(x)} \int_0^x f(y)dy,1 \bigg)_{\omega^{(\beta - 1,\alpha - \beta - 1)}} }{\ds \bigg (\frac1{K(x)},1 \bigg)_{\omega^{(\beta - 1,\alpha - \beta - 1)}}}.
\end{equation}
\end{theorem}

\textbf{Proof}: 
Note that $A$ defined by \eqref{A} is well defined follows from the fact
that $f_1$ and $f_2 \in L^{2}_{\omega^{(\beta - 1,\alpha - \beta - 1)}}(\mrI)$, Lemma \ref{lem:f1f2}.
Also, using this property,
$f_1$ and $f_2$ can be expanded as 
\begin{equation}\label{f1f2Decomp}
f_i = \sum_{j = 0}^{\infty} f_{i,j} \, G_{j}^{(\beta - 1 , \alpha - \beta -1 )}, \quad i = 1, 2
\end{equation}
with
\begin{equation}\label{f1f2Decomp:e1}
f_{i,j} :=  \frac{\big (f_i, G_j^{(\beta - 1 , \alpha - \beta - 1 )} \big )_{\omega^{(\beta - 1  ,  \alpha - \beta - 1)}}}{ \big \|G_j^{(\beta - 1,\alpha-\beta - 1)} \big \|_{\omega^{(\beta - 1  ,  \alpha - \beta - 1)}}^{2}}.      
\end{equation}

Let
\begin{equation}\label{uDecomp}
u(x) := \omega^{(\alpha - \beta  , \beta)}(x) \, \sum_{j = 0}^{\infty} c_j \, G_j^{(\alpha - \beta  , \beta)}(x).
\end{equation}

We use Lemma \ref{lem:beta} to obtain 
\begin{equation}\label{compare1} 
L_r^\alpha u(x) = \sum_{j = 0}^{\infty} c_j L_r^\alpha\big(\omega^{(\alpha - \beta  , \beta)}(x)  G_j^{(\alpha - \beta,\beta)}(x)\big)=-\sum_{j = 0}^{\infty} c_j\lambda_jG_{j+1}^{(\beta-1,\alpha - \beta -1)}(x).
\end{equation}
For $u$ to be a solution to (\ref{EquivForm}), we enforce using (\ref{f1f2Decomp}) 
\begin{equation}\label{compare}\begin{array}{rl}
\ds L_r^\alpha u(x) &\ds \hspace{-0.1in} =-\sum_{j = 0}^{\infty} c_j\lambda_jG_{j+1}^{(\beta-1,\alpha - \beta -1)}(x) = f_1(x) - A f_2(x)\\
&\ds \hspace{-0.1in} = -\sum_{j = 0}^{\infty}( f_{1,j}-Af_{2,j}) \, G_{j}^{(\beta - 1 , \alpha - \beta -1 )}(x) \\
&\ds \hspace{-0.1in} =-(f_{1,0}-Af_{2,0})-\sum_{j = 0}^{\infty}( f_{1,j+1}-Af_{2,j+1}) \, G_{j+1}^{(\beta - 1 , \alpha - \beta -1 )}(x) \ 
\ \mbox{ (using $G_{0}^{(\beta - 1 \, , \, \alpha - \beta - 1)}(x) = 1$)}.
\end{array}\end{equation}
Therefore, for $u$ to be a solution to (\ref{EquivForm}), we require  
\begin{equation}\label{ci}
c_j := (- f_{1, j+1} \ + \ A \, f_{2, j+1})/\lambda_j,~~j\geq 0 
\end{equation}
and 
\begin{equation}\label{ci1}
f_{1,0}-Af_{2,0}= 0. 
\end{equation}  
The value for $A$ given by (\ref{A}) immediately follows from \eqref{ci1}.

To prove that $u(x)$ defined in (\ref{uDecomp}) is a solution of model (\ref{EquivForm})
we define the $N + 1$-term truncation of $u(x)$, $u_N(x)$, as
\begin{equation}\label{u_N}
u_N(x):= \omega^{(\alpha-\beta,\beta)}(x)\sum_{i = 0}^{N } c_{i} \, G_{i}^{(\alpha - \beta \, , \, \beta)}(x).
\end{equation}
Similarly we define
\begin{equation}\label{f_N}
f^N_{k}(x): =  \sum_{i = 0}^{N} f_{k,i} 
G_{i}^{(\beta - 1 \, , \, \alpha - \beta - 1)}(x) ,~~k=1,2. \end{equation}
Now,
\begin{align}
& \| u_M-u_N \|_{ \omega^{(-(\alpha - \beta)  ,   -\beta)}}^{2} \ 
= \ \sum_{i = N+1}^{M} c_{i}^{2} \,  |\| G_{i}^{(\alpha - \beta  ,   \beta)} |\|^{2}  \nonumber  \\
& \quad \le \
\sum_{i = N+1}^{M} \frac{2 f_{1, i+1}^{2}}{\lambda_{i}^{2}} \,  |\| G_{i}^{(\alpha - \beta  ,   \beta)} |\|^{2} 
\ + \
A^{2} \, \sum_{i = N+1}^{M} \frac{2 f_{2, i+1}^{2}}{\lambda_{i}^{2}} \,  |\| G_{i}^{(\alpha - \beta  ,   \beta)} |\|^{2} \nonumber  \\
& \quad \le \
\frac{2}{\lambda_{0}^{2}} \, \sum_{i = N+1}^{M} f_{1, i+1}^{2} \,  |\| G_{i+1}^{(\beta - 1 \,  , \, \alpha -  \beta - 1)} |\|^{2} 
\ + \
\frac{2}{\lambda_{0}^{2}} \, A^{2} \, 
\sum_{i = N+1}^{M} f_{2, i+1}^{2} \,  |\| G_{i+1}^{(\beta - 1 \,  , \, \alpha -  \beta - 1)} |\|^{2} 
\ \ \mbox{(using \eqref{GnNormRelation})} \nonumber 
 \\
& \quad \le \
\frac{2}{\lambda_{0}^{2}} \, \| f^M_{1}-f^N_1 \|_{ \omega^{(\beta - 1 \,  , \, \alpha -  \beta - 1)}}^{2} \ + \
\frac{2}{\lambda_{0}^{2}} \, A^{2} \,  \|f^M_{2}-f^N_2  \|_{ \omega^{(\beta - 1 \,  , \, \alpha -  \beta - 1)}}^{2}  \label{nmes3} 
\end{align}
By Lemma \ref{f1f2} and the assumption that $f\in L^2_{\omega^{(\beta,\alpha-\beta)}}(\mrI)$, 
$\{f_{1}^{N}(x)\}_{N=1}^\infty$ and $\{f_{2}^{N}(x)\}_{N=1}^\infty$ are Cauchy sequences, hence
$\{u_N(x)\}_{N=1}^\infty$ is a Cauchy sequence in $L^2_{\omega^{(-(\alpha-\beta), -\beta)}}(\mrI)$. 
As $L^2_{\omega^{(-(\alpha-\beta), -\beta)}}(\mrI)$ is complete \cite{hes071}, 
$u(x)=\lim_{N\rightarrow \infty}u_N(x)\in L^2_{\omega^{(-(\alpha-\beta), -\beta)}}(\mrI)$. Furthermore, 
\begin{align*}
&\|f_1-Af_2-L_r^\alpha u_N\|_{\omega^{(\beta - 1 \,  , \, \alpha -  \beta - 1)}}\\
&~~~~=\bigg\|f_1-Af_2-L_r^\alpha\bigg( \omega^{(\alpha-\beta,\beta)}\sum_{i = 0}^{N } c_{i} \, G_{i}^{(\alpha - \beta \, , \, \beta)}\bigg)\bigg\|_{\omega^{(\beta - 1 \,  , \, \alpha -  \beta - 1)}}\\
&~~~~=\bigg\|f_1-Af_2- \sum_{i = 0}^{N } (f_{1,i+1}-Af_{2,i+1}) \, G_{i+1}^{(\beta-1,\alpha - \beta -1)}\bigg\|_{\omega^{(\beta - 1 \,  , \, \alpha -  \beta - 1)}}\\
&~~~~=\bigg\|f_1-Af_2- \sum_{i = 0}^{N+1 } (f_{1,i}-Af_{2,i}) \, G_{i}^{(\beta-1,\alpha - \beta -1)}\bigg\|_{\omega^{(\beta - 1 \,  , \, \alpha -  \beta - 1)}} \ \ \mbox{ (using \eqref{ci1})} \\
&~~~~=\|f_1-f_1^{N+1}-A(f_2-f_2^{N+1})\|_{\omega^{(\beta - 1 \,  , \, \alpha -  \beta - 1)}}\\
&~~~~\leq \|f_1-f_1^{N+1}\|_{\omega^{(\beta - 1 \,  , \, \alpha -  \beta - 1)}}
+|A|\|f_2-f_2^{N+1}\|_{\omega^{(\beta - 1 \,  , \, \alpha -  \beta - 1)}} \, .
\end{align*}
As $f_1^{N+1} \stackrel{N \rightarrow \infty}{\longrightarrow} f_{1}$ and 
$f_2^{N+1} \stackrel{N \rightarrow \infty}{\longrightarrow} f_{2}$, then it follows that $L_r^\alpha u \ = \ f_1-Af_2$. \\
\mbox{ } \hfill \qed

\begin{corollary} \label{corunq1}
For $f \in L^{2}_{\omega^{(\beta , \alpha - \beta)}}(0,1)$ and $0 < K_{m} \le K(x) \le K_{M}$ there exists a
unique solution $u \in L^{2}_{\omega^{(-(\alpha - \beta) ,  -\beta)}}(\mrI)$ to \eqref{ModelV}.
\end{corollary}
\textbf{Proof}: From  Theorem \ref{thm:Wellpose} and Lemma \ref{lem:EquivForm} we have the existence of a solution
 $u \in L^{2}_{\omega^{(-(\alpha - \beta) ,  -\beta)}}(\mrI)$ to \eqref{ModelV}. To establish uniqueness, suppose that
 $u_{1}$ and $u_{2} \in L^{2}_{\omega^{(-(\alpha - \beta) ,  -\beta)}}(\mrI)$ are solutions of \eqref{ModelV}. Let 
 $z(x) \ = \ \left( u_{1}(x) \, - \, u_{2}(x) \right) \in L^{2}_{\omega^{(-(\alpha - \beta) ,  -\beta)}}(\mrI)$. Then, $z(x)$ satisfies
\begin{equation}
\begin{array}{c}
\ds L_r^{\alpha}z(x) \ = \ -A/K(x) \, ,~~x \in \mrI \, , \ \ \mbox{ subject to } \ z(0) = z(1) = 0 \, \\[0.1in]
\ds A  = K(0)\big[\big(r{}_0I_x^{2-\alpha} + (1-r){}_xI_1^{2-\alpha}\big)Dz(x)\big]\big|_{x=0}.
\end{array}
\end{equation}
Proceeding in an analogous manner as in the proof of Theorem \ref{thm:Wellpose} we obtain that $z(x) = 0$. Hence 
$u_{1} = u_{2}$.  \\
\mbox{ } \hfill \qed

\subsection{The stability and regularity estimate for the solution to (\ref{EquivForm})}
\begin{theorem}
For $f \in L^{2}_{\omega^{(\beta , \alpha - \beta)}}(0,1)$ and $0 < K_{m} \le K(x) \le K_{M}$, the solution $u$ of (\ref{EquivForm}) has the following stability estimate 
\begin{equation}\label{stability}
\| u \|_{ \omega^{(-(\alpha - \beta)  ,   -\beta)}}\leq C \|f \|_{ \omega^{( \beta,\alpha - \beta)}} \, 
\mbox{  for some constant } C > 0.
\end{equation}
\end{theorem}
\textbf{Proof}:
Note that from \eqref{A}
\begin{equation}
| A | \ \le \ C \, \| f_{1} \|_{ \omega^{(\beta -1 \,  ,  \alpha - \beta - 1)}} \ \le \ C \, \| f \|_{ \omega^{(\beta  ,  \alpha - \beta)}} \, ,
\ \ \mbox{using \eqref{f1f2}}.
\label{bdA}
\end{equation}

From \eqref{nmes3}, with $u_{N} = 0, \, f_{1}^{N} = 0, f_{2}^{N} = 0$ and taking the limit as $M \rightarrow \infty$ we obtain
\begin{align*}
 \| u \|_{ \omega^{(-(\alpha - \beta)  ,   -\beta)}}^{2} &\le  \frac{2}{\lambda_{0}^{2}}
  \left( \| f_{1} \|_{ \omega^{(\beta -1 \,  ,  \alpha - \beta - 1)}}^{2} \
   + \ A^{2} \, \| f_{2} \|_{ \omega^{(\beta -1 \,  ,  \alpha - \beta - 1)}}^{2}  \right) \nonumber \\
&\le \ C \, \|f \|_{ \omega^{( \beta,\alpha - \beta)}}^{2} \, ,
\end{align*}
where, in the last step we have used \eqref{f1f2} and \eqref{bdA}.  \\
\mbox{ } \hfill \qed

The solution of a fractional differential equation, in general, exhibits a lack of regularity at the boundary of the domain. This is
reflected by the weight function $\omega^{(\alpha - \beta , \beta)}(x)$ in the representation for $u(x)$ in \eqref{uDecomp}.
Acknowledging this singular behavior, in the next result we explore the regulaity of $u(x)$ beyond this singular term by 
considering the regularity of $u(x) / \omega^{(\alpha - \beta , \beta)}(x)$.
\begin{theorem} \label{thmrgqd}
For $j \in \mathbb{N}$, if $D^{j}f_k \in L^{2}_{\omega^{(\beta + j - 1 \, , \, \alpha - \beta + j - 1)}}(\mrI)$ for $k=1,2$, then
$D^{j}\big(u(x)/\omega^{(\alpha-\beta,\beta)}(x)\big) \in L^{2}_{\omega^{(\alpha - \beta + j \, , \, \beta + j )}}(\mrI)$.
\end{theorem}
\textbf{Proof}:
Using (\ref{u_N}), (\ref{f_N}) and \eqref{eqC4} we obtain
 \begin{equation}\label{eqC5}  
D^{j} \, f_k^{N}(x) = \ \sum_{i = 0}^{N} f_{k,i}
\frac{\Gamma(i + j + \alpha-1)}{  \Gamma(i + \alpha-1)}
  G_{i - j }^{(\beta + j - 1 \, , \, \alpha - \beta + j - 1)}(x),~~k=1,2,
\end{equation}  
and
\begin{equation}\label{DjuN}
D^{j} \bigg(\frac{u_{N}(x)}{\omega^{(\alpha-\beta,\beta)}}\bigg)
= \ D^{j} \bigg(  \sum_{i = 0}^{N } c_{i} \, G_{i}^{(\alpha - \beta \, , \, \beta)}(x) \bigg)=  \sum_{i = 0}^{N  } c_{i } \, 
\frac{\Gamma(i+j+\alpha+1)}{\Gamma(i+\alpha+1)} \, G_{i-j}^{(\alpha - \beta + j \, , \, \beta + j)}(x) ,  
\end{equation}
where $G_{i}^{(a , b)}(x) = 0$ for $i < 0$. 

Using (\ref{DjuN}) and (\ref{eqC5}),
\begin{align}
&
\bigg\|D^j\bigg(\frac{u_{M}(x)}{\omega^{(\alpha-\beta,\beta)}} \bigg)-D^j\bigg( \frac{u_{N}(x)}{\omega^{(\alpha-\beta,\beta)}} \bigg)\bigg\|^{2}_{\omega^{(\alpha - \beta + j \, , \, \beta + j )}}
 \nonumber  \\
& ~~ = \int_0^1 \omega^{(\alpha - \beta + j \, , \, \beta +j)} \, \bigg(\sum_{i = N+1}^{M} c_{i} \, 
\frac{\Gamma(i+j+\alpha+1)}{\Gamma(i+\alpha+1)} \, G_{i-j}^{(\alpha - \beta +j \, , \, \beta + j)} \bigg)^2  dx  \nonumber \\
&~~=  
 \sum_{i =N+1}^{M} c_{i}^{2} \, 
\bigg( \frac{\Gamma(i+j+\alpha+1)}{\Gamma(i+\alpha+1)} \bigg)^{2} \, |\| G_{i-j}^{(\alpha - \beta + j \, , \, \beta + j)}(x) |\|^{2}  \nonumber \\
&~~\leq C
 \sum_{i =N+1}^{M} \frac{f_{1,i+1}^2+f_{2,i+1}^2}{\lambda_i^2} \, 
\bigg( \frac{\Gamma(i+j+\alpha+1)}{\Gamma(i+\alpha+1)} \bigg)^{2} \, |\| G_{i-j}^{(\alpha - \beta + j \, , \, \beta + j)}(x) |\|^{2} ~ \mbox{(using (\ref{ci}))}\nonumber \\
&~~\leq \frac{C}{\lambda_{0}^2}
 \sum_{i =N+1}^{M} (f_{1,i+1}^2+f_{2,i+1}^2)
\bigg( \frac{\Gamma(i+j+\alpha+1)}{\Gamma(i+\alpha+1)} \bigg)^{2} \, 
\bigg( \frac{(i-j+1)}{(i-j+\alpha)} \bigg)^{2} \, 
|\| G_{i+1-j}^{(\beta + j-1,\alpha - \beta + j-1 )} |\|^{2} ~ \mbox{(using (\ref{GnNormRelation}))}\nonumber \\
&~~\leq C \, 
 \sum_{i =N+1}^{M} (f_{1,i+1}^2+f_{2,i+1}^2)
\bigg( \frac{\Gamma(i+j+\alpha)}{\Gamma(i+\alpha)} \bigg)^{2} \, 
|\| G_{i+1-j}^{(\beta + j-1,\alpha - \beta + j-1 )}(x) |\|^{2} \nonumber \\
&~~=C \bigg( \| D^{j} f_1^{M+1} \, - \, D^{j} f_1^{N+1} \|_{\omega^{(\beta + j - 1 \, , \, \alpha - \beta + j - 1)}}^{2}  
\ + \  \| D^{j} f_2^{M+1} \, - \, D^{j} f_2^{N+1} \|_{\omega^{(\beta + j - 1 \, , \, \alpha - \beta + j - 1)}}^{2}  \bigg) .  \label{thy1} 
\end{align}

As $D^{j}f_k \in L^{2}_{\omega^{(\beta + j - 1 \, , \, \alpha - \beta + j - 1)}}(\mrI)$ for $k=1,2$, 
then $\{D^{j} f^n_{k} \}_{n=1}^\infty$ is a Cauchy sequence in 
$L^{2}_{\omega^{(\beta + j - 1 \, , \, \alpha - \beta + j - 1)}}(\mrI)$, for $k=1,2$. 
Since $\lim_{n\rightarrow\infty}D^j(u_n(x)/\omega^{(\alpha-\beta,\beta)}(x))=D^j(u(x)/\omega^{(\alpha-\beta,\beta)}(x))$, 
then it follows that
$D^{j} \big(u(x)/\omega^{(\alpha-\beta,\beta)}(x)\big ) \in  
L^{2}_{\omega^{(\alpha - \beta + j \, , \, \beta + j )}}(\mrI)$.  \\
\mbox{ } \hfill \qed

\setcounter{equation}{0}
\setcounter{figure}{0}
\setcounter{table}{0}
\setcounter{theorem}{0}
\setcounter{lemma}{0}
\setcounter{corollary}{0}

\section{A spectral approximation method for (\ref{EquivForm})}
In this section, we introduce a spectral approximation method for the model (\ref{EquivForm}) and study its convergence. 
The convergence rate depends on the regularity of $f(x)$ and $K(x)$.

\begin{lemma} \label{hgrreg4f12}
For $f(x) \in H^{s}_{\omega^{(\beta , \alpha-\beta)}}(\mrI)$ and $1/K(x) \in W^{s+1}_\infty(\mrI) $, 
$s \in \mathbb{N}_{0}$,
the following estimates holds
\begin{align}
\| f_{1} \|_{s+1 ,  \omega^{(\beta-1 , \alpha-\beta-1)}}  &\le \ C(s) \, \|1/K\|_{W^{s+1}_\infty}\| f \|_{s,  \omega^{(\beta , \alpha-\beta)}},
\label{bd24f12} \\
\| f_{2} \|_{s+1 ,  \omega^{(\beta-1 , \alpha-\beta-1)}}  &\le \ C(s) \, \|1/K\|_{W^{s+1}_\infty} \, , \label{bdf2}
\end{align}
for constant $C > 0$.
\end{lemma}
\textbf{Proof}: \\
We have that
\begin{align}
\|f_1\|^2_{s+1 ,  \omega^{(\beta-1 , \alpha-\beta-1)}} &= \ 
\sum_{r=0}^{s+1}\bigg\|D^r\bigg(\frac{1}{K(x)}\int_0^xf(y)dy\bigg)\bigg\|^2_{\omega^{(\beta-1+r , \alpha-\beta-1+r)}} \nonumber \\
&= \ \sum_{r=0}^{s+1}\bigg\|\sum_{k=0}^{r}\bigg( \hspace{-0.05in}
\begin{array}{c}
r\\
k
\end{array} \hspace{-0.05in}\bigg)
D^{r-k}\bigg(\frac{1}{K(x)}\bigg)D^k\int_0^xf(y)dy\bigg\|^2_{\omega^{(\beta-1+r , \alpha-\beta-1+r)}} \nonumber \\
&\leq 2\sum_{r=0}^{s+1} \bigg\|D^{r}\bigg(\frac{1}{K(x)}\bigg)\int_0^xf(y)dy\bigg\|^2_{\omega^{(\beta-1+r , \alpha-\beta-1+r)}} \nonumber \\
& \quad \quad + \ 2 \sum_{r=1}^{s+1}\bigg\|\sum_{k=1}^{r}\bigg(\hspace{-0.05in}
\begin{array}{c}
r\\
k
\end{array}\hspace{-0.05in}\bigg)
D^{r-k}\bigg(\frac{1}{K(x)}\bigg)D^k\int_0^xf(y)dy\bigg\|^2_{\omega^{(\beta-1+r , \alpha-\beta-1+r)}} \label{f1:e1} 
\end{align}
For the first term in \eqref{f1:e1}
\begin{align}
& \bigg\|D^{r}\bigg(\frac{1}{K(x)}\bigg)\int_0^xf(y)dy\bigg\|^2_{\omega^{(\beta-1+r , \alpha-\beta-1+r)}}
\ = \ \int_{0}^{1} \, \omega^{(\beta-1+r , \alpha-\beta-1+r)}(x) \, \bigg( D^{r}\bigg(\frac{1}{K(x)}\bigg) \bigg)^{2} \, 
\bigg( \int_0^xf(y)dy \bigg)^{2} \, dx  \nonumber \\
& \quad \quad \le \bigg\| \frac{1}{K(x)} \bigg\|_{W^{r}_{\infty}}^{2} \, 
 \int_{0}^{1} \, \omega^{(\beta-1 , \alpha-\beta-1)}(x) \, 
\bigg( \int_0^xf(y)dy \bigg)^{2} \, dx \ \ \nonumber \\
& \quad \quad \quad \quad  \quad \quad  \mbox{(using  $\omega^{(\beta-1+r , \alpha-\beta-1+r)}(x) \, \le \, 
 \omega^{(\beta-1 , \alpha-\beta-1)}(x)$ )} \nonumber \\
& \quad \quad \le C \, \bigg\| \frac{1}{K(x)} \bigg\|_{W^{r}_{\infty}}^{2} \,  \, \| f \|_{\omega^{(\beta , \alpha-\beta)}} \, , \mbox{using \eqref{fbdn2}}.
\label{bdfret2}
\end{align}
For the second term in \eqref{f1:e1}
\begin{align}
& \bigg\| D^{r-k}\bigg(\frac{1}{K(x)}\bigg) D^k\int_0^xf(y)dy\bigg\|^2_{\omega^{(\beta-1+r , \alpha-\beta-1+r)}}  \nonumber \\
& \quad = \
  \int_{0}^{1} \omega^{(\beta-1+r , \alpha-\beta-1+r)}(x) \, \bigg(  D^{r-k}\bigg(\frac{1}{K(x)}\bigg) \bigg)^{2} \, 
  \bigg( D^k\int_0^xf(y)dy \bigg)^{2} \, dx   \nonumber \\
& \quad  \le \
 \bigg\| \frac{1}{K(x)} \bigg\|_{W^{r - k}_{\infty}}^{2} \, 
  \int_{0}^{1}  \,   \omega^{(\beta-1+r , \alpha-\beta-1+r)}(x) \,
  \bigg( D^{k - 1}  \, f(x) \, \bigg)^{2} \, dx   \nonumber \\
& \quad   \le \
 \bigg\| \frac{1}{K(x)} \bigg\|_{W^{r - k}_{\infty}}^{2} \, 
  \int_{0}^{1}  \,  \omega^{(\beta-1+k , \alpha-\beta-1+k)}(x) \,
  \bigg( D^{k - 1}  \, f(x) \, \bigg)^{2} \, dx   \, , \ \ \mbox{(as $k \le r$)}  \nonumber \\
& \quad   \le \
 \bigg\| \frac{1}{K(x)} \bigg\|_{W^{r - k}_{\infty}}^{2} \, 
  \bigg\| D^{k - 1}  \, f   \bigg\|_{\omega^{(\beta-1+k , \alpha-\beta-1+k)}}^{2} \, . \label{bdrew22}
\end{align}  

Combining \eqref{f1:e1}-\eqref{bdrew22}, we obtain
\begin{align}
\|f_1\|^2_{s+1 ,  \omega^{(\beta-1 , \alpha-\beta-1)}} &\le \ 
C \bigg(  \bigg\| \frac{1}{K(x)} \bigg\|_{W^{s + 1}_{\infty}}^{2} \,  \| f \|_{\omega^{(\beta , \alpha-\beta)}}
 \ + \  \bigg\| \frac{1}{K(x)} \bigg\|_{W^{k + 1}_{\infty}}^{2} \, 
 \bigg( \sum_{k = 1}^{s + 1} \bigg\| D^{k - 1}  \, f   \bigg\|_{\omega^{(\beta-1+k , \alpha-\beta-1+k)}}^{2} \bigg) \bigg)  \nonumber \\
&\le \  C \,  \bigg\| \frac{1}{K(x)} \bigg\|_{W^{s + 1}_{\infty}}^{2} \,  \| f \|_{s , \omega^{(\beta , \alpha-\beta)}}^{2} \, . \nonumber 
\end{align}

Estimate \eqref{bdf2} follows from the definition of $\| \cdot \|_{s , \omega^{(\beta , \alpha-\beta)}}$ and the integrability
of $\omega^{(\beta - 1 + r \, ,  \, \alpha - \beta - 1 + r)}$ on $\mrI$ for $0 \le r \le s+1$. \\
\mbox{ } \hfill \qed

\begin{theorem}   \label{thmuspg}
For $f(x) \in H^{s}_{\omega^{(\beta , \alpha-\beta)}}(\mrI)$ and $1/K(x) \in W^{s+1}_{\infty}(\mrI) $,  
$s \in \mathbb{N}_{0}$,
there exists  $C > 0$ such that
	\begin{align}\label{conv}
	 \| u-u_N\|_{\omega^{(-(\alpha - \beta)  , - \beta)}}
	 \leq \frac{C\|1/K\|_{W^{s+1}_\infty}}{(N+1)^{\alpha+s}}\big(\|f\|_{s,\omega^{(\beta , \alpha - \beta)}} \ +  1 \big).
	\end{align}
\end{theorem}
\textbf{Proof}: We have
\begin{align*}
& \| u-u_N\|_{\omega^{(-(\alpha - \beta)  , - \beta)}}^{2}   \\
&\quad = \int_0^1 \omega^{(-(\alpha - \beta)  , - \beta)} \bigg(\omega^{(\alpha - \beta , \beta)}
\sum_{i= N}^{\infty} c_{i} G_{i}^{(\alpha - \beta , \beta)}\bigg)^2dx     \\
&\quad =\sum_{i= N}^{\infty}c_i^{2} |\| G_{i}^{(\alpha - \beta \, , \, \beta)}  |\|^{2}\leq C\sum_{i= N}^{\infty}\frac{f_{1,i+1}^{2}+f_{2,i+1}^{2} }{\lambda_i^2}|\| G_{i}^{(\alpha - \beta \, , \, \beta)}  |\|^{2} \ \ (\mbox{using } \eqref{ci}) \\
&\quad\leq \frac{C}{(N+1)^{2(\alpha-1)}}\sum_{i= N}^{\infty}(f_{1,i+1}^{2}+f_{2,i+1}^{2})|\| G_{i}^{(\alpha - \beta \, , \, \beta)}  |\|^{2}  \ \  (\mbox{using } \lambda_{i} \sim i^{\alpha - 1})    \\
&\quad \leq \frac{C}{(N+1)^{2(\alpha-1)}}\sum_{i= N}^{\infty}(f_{1,i+1}^{2}+f_{2,i+1}^{2})|\| G_{i+1}^{( \beta -1\, , \, \alpha -\beta-1)}  |\|^{2}  \ \  (\mbox{using } \eqref{GnNormRelation})   \\
&\quad =\frac{C}{(N+1)^{2(\alpha-1)}}
\int_0^1\omega^{( \beta-1,\alpha - \beta-1)}
 \bigg( \sum_{i= N}^{\infty} f_{1,i+1} G_{i+1}^{(  \beta-1,\alpha - \beta-1)} \bigg) ^2dx   \\
 &\quad \quad \quad \quad ~~+\frac{C}{(N+1)^{2(\alpha-1)}}
 \int_0^1 \omega^{( \beta-1,\alpha - \beta-1)}
 \bigg( \sum_{i= N}^{\infty} f_{2,i+1} G_{i+1}^{(  \beta-1,\alpha - \beta-1)}\bigg)^2dx\\
&\quad =\frac{C}{(N+1)^{2(\alpha-1)}}  \bigg( \| f_{1}  - P_{N+1}^{(\beta-1,\alpha-\beta-1)} f_{1} \|_{\omega^{(\beta-1 , \alpha - \beta-1)}}^{2}  \\
&\quad \quad \quad \quad~~ + \   \| f_{2}  - P_{N+1}^{(\beta-1,\alpha-\beta-1)} f_{2} \|_{\omega^{(\beta-1 , \alpha - \beta-1)}}^{2} \bigg) \\
&~~\leq \frac{C}{(N+1)^{2(\alpha+s)}}\big(\|f_1\|^2_{s+1,\omega^{(\beta-1 , \alpha - \beta-1)}}
  +\|f_2\|^2_{s+1,\omega^{(\beta-1 , \alpha - \beta-1)}}\big)  \ \  (\mbox{using } (\ref{Approx}))  \\
&~~\leq  \frac{C\|1/K\|^2_{W^{s+1}_\infty}}{(N+1)^{2(\alpha+s)}}\big(\|f\|_{s,\omega^{(\beta , \alpha - \beta)}} \ +  1 \big)^{2}
 \ \  (\mbox{using } \eqref{bd24f12}, \eqref{bdf2}).
\end{align*}
\mbox{ } \hfill \qed

\begin{corollary}  \label{corL2}
		Under the assumption of Theorem \ref{thmuspg}, the following estimate holds 
		\begin{align}\label{convL2}
		\| u-u_N\|_{L^2}\leq \frac{C\|1/K\|_{W^{s+1}_\infty}}{(N+1)^{\alpha+s}}\big(\|f\|_{s,\omega^{(\beta , \alpha - \beta)}} \ +  1 \big).
		\end{align}
	\end{corollary}
\textbf{Proof}:
	Note that $\omega^{(-(\alpha - \beta) ,-\beta)}\geq 1$. Hence
	\begin{align}
	\|u-u_N\|_{L^2}^2=\int_0^1(u(x)-u_N(x))^2dx\leq \int_0^1\omega^{(-(\alpha - \beta) ,-\beta)}(x)(u(x)-u_N(x))^2dx.
	\end{align}
	Thus (\ref{conv}) leads to (\ref{convL2}).
\mbox{ } \hfill \qed

\setcounter{equation}{0}
\setcounter{figure}{0}
\setcounter{table}{0}
\setcounter{theorem}{0}
\setcounter{lemma}{0}
\setcounter{corollary}{0}

\section{Numerical experiments}

In this section, we present two numerical examples to investigate the sharpness of the estimate (\ref{conv}) and (\ref{convL2}). 

To determine the predicted convergence rate, we need to analyze the regularity of $f(x)$ and $K(x)$. 
For the numerical experiments we will take $K(x) \in C^{\infty}(\mrI)$, bounded, and bounded
away from $0$. Subsequently, we need to determine the largest value of $s$ such that 
$\|f\|_{s,\omega^{(\beta , \alpha-\beta)}}<\infty$. For the two examples considered, 
the most singular parts of $f(x)$ are $x^{1-\alpha}$ and $(1-x)^{1-\alpha}$.
So, we need to find the largest value of $s$ such that $|x^{1-\alpha}|_{s,\omega^{(\beta , \alpha-\beta)}}<\infty$ 
and $|(1-x)^{1-\alpha}|_{s,\omega^{(\beta , \alpha-\beta)}}<\infty$. Note that
\begin{align*}
|x^{1-\alpha}|^2_{s,\omega^{(\beta , \alpha-\beta)}}&=
\|D^sx^{1-\alpha}\|^2_{\omega^{(\beta+s , \alpha-\beta+s)}}
\leq C\int_0^1 (1-x)^{\beta+s} x^{\alpha-\beta+s} \, x^{2(1-\alpha-s)} dx\\
&=C\int_0^1 (1-x)^{\beta+s} x^{2-\alpha-\beta-s} \, dx <\infty,\\
\mbox{(largest value for $s$)} \ \ &\Longrightarrow 2-\alpha-\beta-s>-1 \, , \ \Longrightarrow  s<3-\alpha-\beta \, .
\end{align*}
Similarly,
\begin{align*}
|(1-x)^{1-\alpha}|^2_{s,\omega^{(\beta , \alpha-\beta)}}&=\|D^s(1-x)^{1-\alpha}\|^2_{\omega^{(\beta+s , \alpha-\beta+s)}}\\
&\leq C\int_0^1 (1-x)^{\beta+s} x^{\alpha-\beta+s} \, (1-x)^{2(1-\alpha-s)} \, dx\\
&=C\int_0^1 (1-x)^{2-2\alpha-s+\beta} x^{\alpha-\beta+s}\, dx<\infty,\\
\mbox{(largest value for $s$)} \ \ &\Longrightarrow 2-2\alpha-s+\beta>-1 \, , \ \Longrightarrow  s<3-2\alpha+\beta \, .
\end{align*}
So by (\ref{conv}), the predicted convergence of $\| u-u_N\|_{\omega^{(-(\alpha - \beta)  , - \beta)}}$ is
$ \sim N^{\kappa}$, where $\kappa=\alpha+s<\min\{3-\beta,~3-(\alpha-\beta)\}$. 

\textbf{Remark}: Theorem \ref{thmuspg} and Corollary \ref{corL2} assume that $f(x) \in H^{n}_{\omega^{(\beta , \alpha - \beta)}}(\mrI)$
for $n \in \mathbb{N}$. Here we are assuming that these estimates extend to $f(x) \in H^{s}_{\omega^{(\beta , \alpha - \beta)}}(\mrI)$
for $s \in \mathbb{R}^{+}$.

\textbf{Example 1} Let $K(x)=1$, $\beta$ determined by (\ref{beta}) and
\begin{align*}
 f(x)=\frac{-r}{\Gamma(2-\alpha)}x^{1-\alpha}+\frac{(1-r)}{\Gamma(2-\alpha)}(1-x)^{1-\alpha}.
\end{align*}
The corresponding solution is
\begin{align*}
u(x)=x-\frac{x^{\beta}{}_2 F_1(1+\beta-\alpha,\beta;\beta+1,x)}{{}_2 F_1((1+\beta-\alpha,\beta;\beta+1,1)}.
\end{align*}
 The experimental convergence rate of $\kappa$ for the error in different norms for 
 Example 1 is shown in Table \ref{1,1}, \ref{1,2} and \ref{1,3}.

\begin{table}[H]
	\setlength{\abovecaptionskip}{0pt}
	\centering
	\caption{Example 1 with $\alpha=1.70$, $r=0.34$ and $\beta=0.90$.}	\label{1,1}
	\vspace{0.5em}	
	\begin{tabular}{ccccc}
		\hline
		$N$&$\|u-u_N\|_{L^2_{\omega^{(-(\alpha-\beta),-\beta)}}}$ &$\kappa$&$\|u-u_N\|_{L^2}$ &$\kappa$\\
		\cline{1-5}
16&	2.69E-04&		&5.33E-05&	\\
18&	2.12E-04&	2.11 &	4.03E-05&	2.51\\ 
20&	1.72E-04&	2.12 &	3.13E-05&	2.52 \\
22	&1.42E-04&	2.12 &	2.49E-05&	2.53 \\
24	&1.19E-04&	2.12 &	2.01E-05&	2.54 \\
		\hline
		Pred.&&2.10&&2.10\\
		\hline	
	\end{tabular}
\end{table}

\begin{table}[H]
	\setlength{\abovecaptionskip}{0pt}
	\centering
	\caption{Example 1 with $\alpha=1.40$, $r=0.62$ and $\beta=0.60$.}	\label{1,2}
	\vspace{0.5em}	
	\begin{tabular}{ccccc}
		\hline
		$N$&$\|u-u_N\|_{L^2_{\omega^{(-(\alpha-\beta),-\beta)}}}$ &$\kappa$&$\|u-u_N\|_{L^2}$ &$\kappa$\\
		\cline{1-5}
16&	2.89E-04&		&7.11E-05&	\\
18	&2.25E-04&	2.25 &	5.30E-05&	2.63\\ 
20	&1.80E-04&	2.25 &	4.07E-05&	2.64 \\
22	&1.46E-04&	2.25 &	3.20E-05&	2.65 \\
24	&1.21E-04&	2.25 &	2.56E-05&	2.66 \\
	\hline
		Pred.&&2.20&&2.20\\
		\hline	
	\end{tabular}
\end{table}

\begin{table}[H]
	\setlength{\abovecaptionskip}{0pt}
	\centering
	\caption{Example 1 with $\alpha=1.70$, $r=0.50$ and $\beta=0.85$.}	\label{1,3}
	\vspace{0.5em}	
	\begin{tabular}{ccccc}
		\hline
		$N$&$\|u-u_N\|_{L^2_{\omega^{(-(\alpha-\beta),-\beta)}}}$ &$\kappa$&$\|u-u_N\|_{L^2}$ &$\kappa$\\
		\cline{1-5}
16&	2.80E-04&		&5.46E-05&	\\
18	&2.22E-04&	2.10 &	4.13E-05&	2.50\\ 
20&	1.80E-04	&2.11 &	3.21E-05&	2.51 \\
22	&1.48E-04	&2.11 	&2.55E-05&	2.52 \\
24	&1.24E-04	&2.11 	&2.07E-05&	2.53 \\
		\hline
		Pred.&&2.15&&2.15\\
		\hline	
	\end{tabular}
\end{table}

\textbf{Example 2} Let $K(x)=e^x$ and
\begin{align*}
 f(x)&=-re^x\left(\frac{x^{1-\alpha}}{\Gamma(2-\alpha)}-\frac{x^{2-\alpha}}{\Gamma(3-\alpha)}-2\frac{x^{3-\alpha}}{\Gamma(4-\alpha)}\right)\\
 &-(1-r)e^x\left(\frac{(1-x)^{1-\alpha}}{\Gamma(2-\alpha)}-3\frac{(1-x)^{2-\alpha}}{\Gamma(3-\alpha)}+2\frac{(1-x)^{3-\alpha}}{\Gamma(4-\alpha)}\right).
\end{align*}
The corresponding solution is $u=x(1-x)$.
The experimental convergence rate of $\kappa$ for the error in different norms for 
Example 2 is shown in Table \ref{2,1}, \ref{2,2} and \ref{2,3}.

\begin{table}[H]
 	\setlength{\abovecaptionskip}{0pt}
 	\centering
 	\caption{Example 2 with $\alpha=1.70$, $r=0.34$ and $\beta=0.90$.}	\label{2,1}
 	\vspace{0.5em}	
 	\begin{tabular}{ccccc}
 		\hline
 		$N$&$\|u-u_N\|_{L^2_{\omega^{(-(\alpha-\beta),-\beta)}}}$ &$\kappa$&$\|u-u_N\|_{L^2}$ &$\kappa$\\
 		\cline{1-5}
16&	2.40E-04&		&4.69E-05&	\\
18	&1.92E-04&	2.02 &	3.59E-05&	2.41\\ 
20	&1.56E-04&	2.03 &	2.82E-05&	2.43 \\
22	&1.30E-04&	2.04 &	2.25E-05&	2.45 \\
24	&1.10E-04&	2.05 &	1.84E-05&	2.46 \\
 	\hline
 	Pred.&&2.10&&2.10\\
 		\hline	
 	\end{tabular}
 \end{table}
 
\begin{table}[H]
 	\setlength{\abovecaptionskip}{0pt}
 	\centering
 	\caption{Example 2 with $\alpha=1.40$, $r=0.62$ and $\beta=0.60$.}	\label{2,2}
 	\vspace{0.5em}	
 	\begin{tabular}{ccccc}
 		\hline
 		$N$&$\|u-u_N\|_{L^2_{\omega^{(-(\alpha-\beta),-\beta)}}}$ &$\kappa$&$\|u-u_N\|_{L^2}$ &$\kappa$\\
 		\cline{1-5}
 16&	2.55E-04&	&	6.13E-05&\\	
18&	2.01E-04&	2.14 &	4.63E-05&	2.51\\ 
20&	1.62E-04&	2.15 &	3.60E-05&	2.53 \\
22&	1.33E-04&	2.16 &	2.85E-05&	2.55 \\
24&	1.11E-04&	2.17 &	2.30E-05&	2.56 \\
 		\hline
 		Pred.&&2.20&&2.20\\
 		\hline	
 	\end{tabular}
 \end{table}
 
 \begin{table}[H]
 	\setlength{\abovecaptionskip}{0pt}
 	\centering
 	\caption{Example 2 with $\alpha=1.40$, $r=0.50$ and $\beta=0.7$.}	\label{2,3}
 	\vspace{0.5em}	
 	\begin{tabular}{ccccc}
 		\hline
 		$N$&$\|u-u_N\|_{L^2_{\omega^{(-(\alpha-\beta),-\beta)}}}$ &$\kappa$&$\|u-u_N\|_{L^2}$ &$\kappa$\\
 		\cline{1-5}
16&	2.73E-04&		&6.54E-05&	\\
18&	2.14E-04&	2.16 &	4.95E-05&	2.51\\ 
20&	1.73E-04&	2.17 &	3.84E-05&	2.53 \\
22&	1.42E-04&	2.18 &	3.04E-05&	2.55 \\
24&	1.18E-04&	2.19 &	2.46E-05&	2.57 \\
 		\hline
 		Pred.&&2.30&&2.30\\
 		\hline	
 	\end{tabular}
 \end{table}
 
 \setcounter{equation}{0}
\setcounter{figure}{0}
\setcounter{table}{0}
\setcounter{theorem}{0}
\setcounter{lemma}{0}
\setcounter{corollary}{0}

\section{Conclusion}
\label{secConc}
In this paper we have established the wellposedness of the variable coefficient fractional differential equation \eqref{ModelV}
for $f \in L^{2}_{\omega^{(\beta , \alpha - \beta)}}(\mrI)$. A spectral approximation method was proposed and analyzed. For the error
between the approximation and the true solution, measured in the appropriately weighted $L^{2}$ norm, the numerical results are
in very good agreement with the derived theoretical estimates. Additionally, the numerical results indicate that the established
convergence rate for the error in the $L^{2}$ norm (Corollary \ref{corL2}) may not be optimal.

 \section*{Acknowledgements}
 This work was funded by the OSD/ARO MURI Grant W911NF-15-1-0562 and by the National Science Foundation under Grant DMS-1620194.

\end{document}